\documentclass[preprint,12pt]{elsarticle}
\usepackage[top=1.3in, bottom=2cm, left=2.5cm, right=2.5cm]{geometry}

\usepackage{graphicx,epsfig}

\usepackage{amssymb,amsmath}
\usepackage{graphicx}

\journal{}
\begin{document}
\begin{frontmatter}
\newtheorem{thm}{Theorem}[section]
\newtheorem{lem}{Lemma}[section]
\newtheorem{coro}{Corollary}[section]
\newtheorem{prop}{Proposition}[section]
\newdefinition{definition}{Definition}
\newdefinition{rmk}{Remark}
\newproof{pf}{Proof}
\newtheorem{ex}{Example}
\title{An Unknotting Sequence for Torus Knots\tnoteref{t1}}
\tnotetext[t1]{This document is a collaborative effort.}

\author[v]{Vikash. S}

\author[mp]{Madeti. P\corref{cor1}}
\ead{prabhakar@iitrpr.ac.in}
\ead[url]{http://www.iitrpr.ac.in/html/faculty/prabhakar.shtml}
\cortext[cor1]{Corresponding author}
\address[v]{Department of Mathematics, IIT Ropar, Rupnagar- 140001, India.}
\address[mp]{Room No. 207, Department of Mathematics, IIT Ropar, Rupnagar - 140001, India.}

\begin{abstract}
In this paper, the authors give an unknotting sequence for torus knots and also determine the unknotting numbers of $_n14_{17191}, \ _n14_{14274},
\ _n14_{18351}, \ _n14_{24498}$ and some other knots from the knot table of Hoste-Thistlethwite.
\end{abstract}

\end{frontmatter}

\section{Introduction}
An unknotting sequence for a knot or a link $K$ is a finite sequence of knots or links 
\[K=K_n, K_{n-1}, K_{n-2},\cdots, K_1, K_0=\ trivial\ link,\]
such that:
\begin{enumerate}
 \item The unknotting number of $K_i$ is $i$, i.e. $u(K_i)=i,\ 0\leq i \leq n$, and
 \item two succeeding knots or links of the sequence are related by one crossing change.
\end{enumerate}

Even though every knot in $S^3$ can be unknotted by a finite sequence of crossing changes, it is interesting to see that every knot of 
unknotting number at least two can be unknotted via infinitely many different knots of unknotting number one \cite{noteoncro}. Sebastian, in his 
paper \cite{unknottingseq}, showed that the unknotting number of a quasipositive knot is equal to its genus if and only if it lies in an unknotting 
sequence of some torus knot. However, many non-quasipositive knots also exist in unknotting sequences of torus knots. For example, the 
non-quasipositive knots $8_2$ and $8_7$ exist in unknotting sequences $5_1, 8_2, 0_1$ and $5_1, 8_7, 0_1$ 
respectively.

In \cite{unknottingpro}, the authors presented a new approach to unknot torus knots and extended the same to torus links. In particular, the authors introduced unknotting crossing data and minimal unknotting crossing data in \cite{unknottingpro}. This minimal unknotting crossing data helps in selecting a pattern of crossings from a toric braid representation of torus knots, such that switching of all crossings at this selected crossing data results a braid whose closure is isotopically equivalent to the trivial knot. Based on this selection of pattern of crossings, authors determine the exact braid representation of all intermediate braids, whose closures give an unknotting sequence for torus knots.

 In Section \ref{sec2}, the authors give an unknotting sequence for torus knots. In Section \ref{sec3}, the unknotting numbers of $_n13_{604}, \ _n14_{17191}, \ _n14_{14274},\ _n14_{18351}, \ _n14_{24498}$ and some other knots is obtained by showing that each of these knots lies in some unknotting sequence of torus knots.
 
In Subsection \ref{subsec3.1}, a sharp upper bound for the unknotting number of two special classes of knots has been discussed.

\section{An Unknotting Sequence of Torus Knots}\label{sec2}
The unknotting number of torus knots is well known \cite{km-2, km-3}. Here we give an unknotting sequence for torus knots. Throughout this paper, we use term torus knots for both torus knots and torus links. We consider torus knots as the closure of $(\sigma_1\sigma_2\cdots \sigma_{p-1})^q$ and denote as  $K(p,q) = cl(\sigma_1\sigma_2\cdots \sigma_{p-1})^q$. We denote $1+2+\cdots +n$ as $\sum n$ and use two braid types\\ $\eta_l =\left\{
\begin{array}{c l l}      
    \sigma_{k+2-l}\sigma_{k+3-l}\cdots\sigma_{p-l} & if\ 1\leq l < j\\
    \sigma_{k+3-l}\sigma_{k+4-l}\cdots\sigma_{p-l} & if\ j< l < k+2
   \end{array}\right.$ and\ 
  $\beta_j = \sigma_{k+3-j}\sigma_{k+4-j}\cdots\sigma_{p-1},$\\
  which depends on $j$ and $k$ (we will define $j$ and $k$ whenever we use $\eta_l$ and $\beta_j$). 

To find an unknotting sequence of torus knots, we divide all torus knots in two classes
\begin{enumerate}
\item when $q\equiv 0\ or\ \pm1\  (mod\ p)$
\item otherwise
\end{enumerate}
\begin{rmk}\label{remark1}
Observe that for a sequence of knots $K_{n}, K_{n-1},\cdots K_{m},\cdots, K_{0}$,  where $K_0$ is trivial knot and $K_i$ can be obtained from $K_{i+1}$ by one crossing change for any $i$, if $ u(K_m)=m $ for some $m$ then $K_{m}, K_{m-1},\cdots, K_{0}$ is an unknotting sequence for $K_{m}$.
\end{rmk}

\begin{thm}\label{thm1}
Let $q = p\ or\ p\pm1$, $n=u(K(p,q))$ and for any $i\leq n$
\[ K_{n-i}=cl(\eta_1\eta_2\cdots\eta_{j-1}\beta_j\eta_{j+1}\cdots\eta_{k+2}(\sigma_1\sigma_2\cdots \sigma_{p-1})^{q-(k+2)})\]
where, 
\begin{align}
k &=sup\{n\in \mathbb{Z^+} : \sum n <i\}, \notag \\
j &= k+2-(i-\sum k), \notag 
\end{align}
then
\[K_n, K_{n-1}, K_{n-2},\cdots, K_1, K_0=trivial\ knot,\]
is an unknotting sequence for $K(p,q)$ torus knot.
\end{thm}
\begin{pf}
Observe that $K_n=K(p,q)$ and $u(K_n)=u(K(p,q))=n$, so by Remark \ref{remark1}, $K_n, K_{n-1}, K_{n-2},\cdots, K_1, K_0$ will be an unknotting sequence of $K(p,q)$ if, for each $i<n$, $K_{n-(i+1)}$ is obtained from $ K_{n-i} $ by making one crossing change.\\
If $i \neq \sum (k+1)$, then
\begin{align}
K_{n-i}&=cl(\eta_1\eta_2\cdots\eta_{j-1}\beta_j\eta_{j+1}\cdots\eta_{k+2}(\sigma_1\sigma_2\cdots \sigma_{p-1})^{q-(k+2)}) \notag \\
&= cl(\eta_1\eta_2\cdots\eta_{j-2}\underbrace{\sigma_{k+3-j}\sigma_{k+4-j}\cdots\sigma_{p+1-j}}\underbrace{\sigma_{k+3-j}\sigma_{k+4-j}\cdots\sigma_{p-1}}\notag \\ &~~~~\eta_{j+1} \cdots\eta_{k+2}(\sigma_1\sigma_2\cdots \sigma_{p-1})^{q-(k+2)}) \notag \\
&= cl(\eta_1\eta_2\cdots\eta_{j-2}\sigma_{k+3-j}\underbrace{\underline{\sigma_{k+3-j}}\sigma_{k+4-j}\cdots\sigma_{p-1}}\underbrace{\sigma_{k+3-j}\sigma_{k+4-j}\cdots\sigma_{p-j}}\eta_{j+1}\notag \\ &~~~~\cdots\eta_{k+2}(\sigma_1\sigma_2\cdots \sigma_{p-1})^{q-(k+2)}). \notag
\end{align}
Then after making crossing change at the underlined position $\underline{\sigma_{k+3-j}}$, we obtain\\
$ cl(\eta_1\eta_2\cdots\eta_{j-2}\underbrace{\sigma_{k+4-j}\cdots\sigma_{p-1}}\underbrace{\sigma_{k+3-j}\sigma_{k+4-j}\cdots\sigma_{p-j}}\eta_{j+1} \cdots\eta_{k+2}(\sigma_1\sigma_2\cdots \sigma_{p-1})^{q-(k+2)})$\\
$= cl(\eta_1\eta_2\cdots\eta_{j-2}\beta_{j-1}\eta_{j}\cdots\eta_{k+2}(\sigma_1\sigma_2\cdots \sigma_{p-1})^{q-(k+2)})
 = K_{n-(i+1)}.$\\

\noindent If $i = \sum (k+1)$, then
 \begin{align}
K_{n-i}&=cl(\beta_1\eta_2\cdots\eta_{k+2}(\sigma_1\sigma_2\cdots \sigma_{p-1})^{q-(k+2)}) \notag \\
&= cl(\beta_1\eta_2\cdots\eta_{k+1}\sigma_{1}\sigma_{2}\cdots\sigma_{p-(k+2)}(\sigma_1\sigma_2\cdots \sigma_{p-1})^{q-(k+2)}) \notag \\
&= cl(\beta_1\eta_2\cdots\eta_{k+1}\sigma_{1}\underline{\sigma_{1}}\sigma_{2}\cdots\sigma_{p-1}\sigma_{1}\sigma_{2}\cdots
\sigma_{p-(k+3)}(\sigma_1\sigma_2\cdots \sigma_{p-1})^{q-(k+3)}). \notag 
\end{align}
Then after making crossing change at the underlined position $\underline{\sigma_{1}}$, we obtain\\
\[cl(\beta_1\eta_2\cdots\eta_{k+1}\sigma_{2}\cdots\sigma_{p-1}\sigma_{1}\sigma_{2}\cdots
\sigma_{p-(k+3)}(\sigma_1\sigma_2\cdots \sigma_{p-1})^{q-(k+3)})\]
\[= cl(\eta_1\eta_2\cdots\eta_{k+1}\beta_{k+2}\eta_{k+3}(\sigma_1\sigma_2\cdots \sigma_{p-1})^{q-(k+3)})=K_{n-(i+1)}.\]
Hence proved.
\end{pf}
\begin{thm}\label{thm2}
Let $q\equiv 0\ or\ \pm1\ (mod\ p)$, $n=u(K(p,q))$ and for $i\leq n$
\[ K_{n-i}=cl(\eta_1\eta_2\cdots\eta_{j-1}\beta_j\eta_{j+1}\cdots\eta_{k+2}(\sigma_1\sigma_2\cdots \sigma_{p-1})^{q-(mp+k+2)})\]
where
\begin{align}
m &= [i/\sum(p-1)]\notag \\
k &= sup\{n\in \mathbb{Z^+} : \sum n <i\ mod\sum(p-1)\}\notag \\
j &= (k+2)-(i-m\sum(p-1)-\sum k),\notag
\end{align}
then 
$K_n, K_{n-1}, K_{n-2},\cdots, K_1, K_0=trivial\ knot$, is an unknotting sequence for $K(p,q)$.
\end{thm}
\begin{pf}
Let $q=mp\ or\ mp\pm 1,$ then based on the proof of Theorem \ref{thm1}, it is easy to observe that after changing $\sum(p-1)$ crossings we have the torus knot $K_{n-\sum(p-1)}=K(p,q-p)$. Now again applying $\sum(p-1)$ crossing changings from $K(p,q-p)$, we have $K(p,q-2p)$. Continuing in this procedure, we will get $K(p,p)$ or $K(p,p\pm 1)$. The result follows by Theorem \ref{thm1}.
\end{pf}

Now we give an unknotting sequence for the remaining torus knots $K(p,q)$, i.e., the case when $q\equiv a\ (mod\ p)$, where $a\neq 0, \pm 1.$\\
By considering the torus knot of type $K(p,rp+a)$ where $0<a<p$, we will show that, after choosing the pattern and making $r\sum (p-1) + \sum (a-1)$ crossing changes accordingly, $K(p-a,a) $ is left.

In Theorem \ref{thm3}, we will show a partial unknotting sequence for the torus knot of type $K(p,rp+a)$, which will end at $K(p-a,a)$. 
\begin{rmk}\label{remark2} 
An unknotting sequence for the torus knot $K(p,rp+a)$ can be obtained by applying the following two steps:\\
(i) Use Theorem \ref{thm3} for finitely many times, until we get a sequence which ends at $K(d,md)$ or $ K(d,md\pm1) $, where $d=gcd(p,rp+a)$.\\
(ii) After completion of step (i), use Theorem \ref{thm2} on $K(d,md)$ or  $ K(d,md\pm 1) $.
\end{rmk}
\begin{thm}\label{thm3} Let $q=rp+a$ for some $r\geq 0$ and $0<a< p$ and for $ i\leq r\sum (p-1) + \sum (a-1)$ $K_{n-i}$ is same as in Theorem~\ref{thm2}, then 
\[K_n, K_{n-1}, K_{n-2},\cdots, K_{n-(r\sum (p-1) + \sum (a-1)-1)}, K_{n-(r\sum (p-1) + \sum (a-1))}\]
is a part of an unknotting sequence for $K(p,q)$. Moreover, $K_{n-(r\sum (p-1) + \sum (a-1))} = K(p-a,a)$.
\end{thm}
\begin{pf}
As proved in Theorem \ref{thm2}, we can obtain $K_{n-i}$ from $K_{n-(i-1)}$ with one crossing change for each $i\leq r\sum (p-1) + \sum (a-1)$. Note that
\begin{align}
K_{n-(r\sum (p-1) + \sum (a-1))} &= cl(\sigma_a\sigma_{a+1}\cdots\sigma_{p-1}\sigma_{a-1}\sigma_{a}\cdots
\sigma_{p-2}\cdots\sigma_{1}\sigma_{2}\cdots\sigma_{p-a}) \notag \\
&= cl(\sigma_{1}\sigma_{2}\cdots\sigma_{p-a-1})^{a}. \notag \\
&= K(p-a,a) \notag
\end{align}
Now to show $ K_n, K_{n-1}, K_{n-2},\cdots, K_{n-(r\sum (p-1) + \sum (a-1)-1)}, K_{n-(r\sum (p-1) + \sum (a-1))} $ is a part of an unknotting sequence of $K(p,q)$, observe that \\
\begin{align}
u(K(p,q))-total\ number\ of\ crossing\ changed \notag \\
= u(K(p,q))-(r\sum (p-1) &+ \sum (a-1)) = \frac{(a-1)(p-a-1)}{2},\notag
\end{align}
which is the unknotting number of $K(p-a,a)$, the last term of the sequence. Thus, for each $i$, the unknotting number of $K_i$ is $i$. Hence the given sequence is a part of an unknotting sequence of $K(p,q)$. 
\end{pf}
\begin{ex}
An Unknotting Sequence for $K(5,7)$:\\
By Remark \ref{remark2}, first we find a part of an unknotting sequence for $K(5,7)$ by making $11$ crossing changes. A partial unknotting sequence for $K(5,7)$ will be \[ K_{12}, K_{11}, K_{10}, K_{9}, K_{8}, K_{7}, K_{6}, K_{5}, K_{4}, K_{3}, K_{2}, K_{1} = K(3,2).\] Here, for $2<i\leq 12$,
\[ K_i = cl(\eta_1\eta_2\cdots\eta_{j-1}\beta_j\eta_{j+1}\cdots\eta_{k+2}(\sigma_1\sigma_2\cdots \sigma_{4})^{5-k})\]
where
\begin{align}
k &= sup\{n\in \mathbb{Z^+} : \sum n <12-i\}\notag \\
j &= \sum k+k +i -10\notag
\end{align}
and
\begin{align}
K_2 &= (\sigma_1\sigma_2\sigma_3\sigma_4)^2 = K(5,2)\notag \\
  K_1 &= \sigma_2\sigma_3\sigma_4\sigma_1\sigma_2\sigma_3 = K(3,2)\notag  
\end{align}
Now applying Theorem \ref{thm2} on $K(3,2)$, an unknotting sequence for $K(5,7)$ is
\[ K_{12}, K_{11}, K_{10}, K_{9}, K_{8}, K_{7}, K_{6}, K_{5}, K_{4}, K_{3}, K_{2}, K_{1}, K_0 = trivial\ knot.\]
\end{ex}

\section{Unknotting Number of Some Knots}\label{sec3}
In this section we determine the unknotting numbers for $_n13_{604},\ _n14_{17191},\ _n14_{14274},\ _n14_{18351}$, $ _n14_{24498} $ and some other knots. In subsection \ref{subsec3.1}, we give a sharp upper bound for the unknotting number of some knot classes using the results from Section \ref{sec2}. Unknotting sequences of $K(4,7)$ can be obtained by making crossing changes corresponding to Minimal unknotting crossing data $\left<6,8,9,10,11,12,18,20,21\right>$. Observe that $_n10_{21}$ lies in an unknotting sequence of $K(4,7)$ and can be obtained by changing crossings $\left<6,9,12,18,20\right>$ from Minimal unknotting crossing data of $K(4,7)$. So $u(_n10_{21})=u(K(4,7)) - \#\left<6,9,12,18,20\right> = 4$. In a similar manner, we determine the unknotting number of some knots which lies in some unknotting sequence of torus knots. After changing a fixed random crossing selection in Minimal unknotting crossing data of the torus knot $K(p,q)$, the unknotting number of the resultant knot is equal to $u(K(p,q))-\# (random\ crossing\ selection)$. Some observations are given in Table~\ref{table1}. For example, the unknotting number of $_n14_{17191}$ is $5$, since this knot is equivalant to closure of a toric braid obtained from $K(4,7)$ after changing $\left<6,9,12,21\right>$ crossings from Minimal unknotting crossing data of $K(4,7)$.
\begin{table}
\begin{center}
\caption{Unknotting number of some knots}
\begin{tabular}{| c| c | c| c | l |}
\hline
Knot & Unknotting & Torus & Minimal Unknotting & Crossing changed\\
H-T Notation  &  Number & Knot & Crossing Data &\\
\hline
$_a9_{38}$ &$3$ &$K(5,6)$  & $ 8,11,12,14,15,16,17,18,19,20 $ & $8,11,12,14,16,19,20$\\
$_n10_{21}$ &$4$ &$K(4,7)$ & $6,8,9,10,11,12,18,20,21$ &$6,9,12,18,20$\\
$_n12_{417}$ &$4$ &$K(4,7)$ & $6,8,9,10,11,12,18,20,21$ &$6,9,12,18,21$\\
$_n13_{604}$ &$5$ &$K(4,7)$ & $6,8,9,10,11,12,18,20,21$ &$6,9,11,21$\\
$_n14_{14274}$ &$4$ &$K(5,6)$ & $ 8,11,12,14,15,16,17,18,19,20 $ &$8,11,12,16,18,20$\\
$_n14_{17191}$ &$5$ &$K(4,7)$ & $6,8,9,10,11,12,18,20,21$ &$6,9,12,21$\\
$_n14_{18351}$ &$4$ &$K(5,6)$ & $ 8,11,12,14,15,16,17,18,19,20 $ &$8,11,12,15,16,20$\\
$_n14_{24498}$ &$5$ &$K(5,6)$ & $ 8,11,12,14,15,16,17,18,19,20 $ &$8,11,12,16,20$\\
\hline
\end{tabular} \label{table1}
\end{center}
\end{table}

\subsection{Sharp Upper Bound for the Unknotting Number of Some Knot Classes} \label{subsec3.1}
Here we consider two knot classes and give sharp upper bound for the unknotting number of knots in these classes. This upper bound is based on the unknotting procedure given in \cite{unknottingpro}.\
\begin{thm}\cite[Theorem 3.1]{unknottingpro}\label{[1]-thm3.1}
 For every $n$, the $(n+1)$-braid
\[\sigma_1\sigma_2\cdots\sigma_n \sigma_1\sigma_2\cdots\sigma_{n-1}\sigma_n^{-1}
\sigma_1\sigma_2\cdots\sigma_{n-1}^{-1}\sigma_n^{-1}\cdots \sigma_1^{-1}\sigma_2^{-1}\cdots\sigma_n^{-1}\]
is a trivial $(n+1)$-braid.
\end{thm}
\begin{thm}\cite[Theorem 3.2]{unknottingpro}\label{[1]-thm3.2}
 Let $K(p, q)$ be a torus knot with $(p, q) = 1$. If \[q\equiv 1 \ \textrm{or} \ p-1 \ (mod \ p),\] then the $U-$crossing data for $B(p, q)$
is a minimal unknotting crossing data for $B(p, q)$ (or $K(p, q)$).
\end{thm}
Note: For  $U-$crossing data and other details see \cite{unknottingpro}.\\

Using these results, we give sharp upper bound for the following two classes of knots:
\begin{enumerate}
\item Consider quasi-toric braid knots with representation $K= cl(\sigma_1\sigma_2^{-1})^q$, where $q=6m+r$; then
\[u(K)\leq \left\{
\begin{array}{c l c}      
    4m & if\ r=0,1\\
    4m+1 & if\ r=2\\
    4m+2 & if\ r=3,4,5
\end{array}\right.\]
This is easy to observe that from $cl(\sigma_1\sigma_2^{-1})^q$, where $q=6m+r$, if we change\\
$
    \left<12(i-1)+ \{2, 5, 8, 11\}\right>\   for\ r\ =\ 0,\ 1;\ 
    \left<12(i-1)+ \{2, 5, 8, 11\}, 12m+2\right>\  for\ r\ =\ 2;
    \left<12(i-1)+ \{2, 5, 8, 11\}, 12m+\{2, 5\}\right>\  for\ r\ =\ 3,\ 4;$ and \\
    $\left<12(i-1)+ \{2, 5, 8, 11\}, 12m+ \{4, 9\} \right>\  for\ r\ =\ 5;
$ where $1\leq i \leq m$, then by Theorem \ref{[1]-thm3.1}, Theorem \ref{[1]-thm3.2} and some braid relations, the closure of the resultant braid is equivalant to unknot.

\hspace*{0.5cm} In case when $q=2,4$ and $5$, the above inequality gives exact unknotting numbers. The corresponding knots are $4_{1}, 8_{18}$ and $10_{123}$ and their unknotting numbers are $1$, $2$ and $2$ respectively. It is interesting to note that these are the only knots in this category whose unknotting numbers are known.\\
For $q=7$, $u(_a14_{19470})\leq 4$ and for $q=8$, $u(_a16_{379778})\leq 5$.

\item Consider quasi-toric braid knots with representation $K=cl(\sigma_1\sigma_2\sigma_3^{-1})^q$, where $q=4m+r$; then
\[u(K)\leq \left\{
\begin{array}{c l}      
    q & if\ 4|q\\
    q-1 & otherwise
\end{array}\right.\]
This is easy to observe that from $cl(\sigma_1\sigma_2\sigma_3^{-1})^q$, where $q=4m+r$, if we change\\
$ 
    \left<12(i-1)+ \{3, 8, 10, 11\}\right>\  for\ r\ =\ 0,\ 1;
    \left<12(i-1)+ \{3, 8, 10, 11\}, 12m+5\right> \  for\ r\ =\ 2;\ and\
    \left<12(i-1)+ \{3, 8, 10, 11\}, 12m+\{3, 8\}\right>\ for\ r\ =\ 3;
$ crossings for $1\leq i \leq m$,  then by Theorem \ref{[1]-thm3.1}, Theorem \ref{[1]-thm3.2} and some braid relations the closure of the resultant braid is equivalant to unknot.

\hspace*{0.5cm} In case when $q=3$, this inequality gives exact unknotting number. The corresponding knot is $9_{47}$ and its unknotting number is $2$.  It is interesting to note that this is the only knot in this category whose unknotting number is known.
For $q=5$
\[u(_n15_{166130})\leq 4.\]
\end{enumerate}

\noindent \textbf{Acknowledgements}\\

Authors thank the reviewer for his valuable comments and suggestions. Also the first author thanks CSIR, New Delhi and IIT Ropar for providing financial assistance and research facilities.\\









\end{document}